\documentclass[10pt]{article}
\usepackage{color}
\usepackage{amssymb}
\usepackage{amsthm,array,amssymb,amscd,amsfonts,latexsym, url}
\usepackage{amsmath}
\usepackage[all]{xy}

\newtheorem{theo}{Theorem}[section]
\newtheorem{prop}[theo]{Proposition}
\newtheorem{lemm}[theo]{Lemma}
\newtheorem{coro}[theo]{Corollary}
\newtheorem{rema}[theo]{Remark}

\newtheorem{conj}[theo]{Conjecture}

\title{ Bloch's conjecture for
Catanese and Barlow surfaces}
\author{Claire Voisin
\\CNRS and \'{E}cole Polytechnique}
\date{}

\newcommand{\cqfd}
{%
\mbox{}%
\nolinebreak%
\hfill%
\rule{2mm}{2mm}%
\medbreak%
\par%
}
\newfont{\gothic}{eufb10}
\begin{document}
\maketitle
\setcounter{section}{-1}
\begin{flushright} \`{A} la m\'{e}moire de Friedrich Hirzebruch
\end{flushright}
\begin{abstract}  Catanese surfaces are regular surfaces of general
type with $p_g=0$. They specialize to double covers of Barlow
surfaces. We prove that the $CH_0$ group of a  Catanese surface is
equal to $\mathbb{Z}$, which implies the same result for the Barlow
surfaces.
\end{abstract}

\section{Introduction}
In this paper, we establish an improved version of  the main theorem
of \cite{genhodgebloch} and use it in order to prove the Bloch
conjecture for Catanese surfaces. We will also give a conditional
application (more precisely, assuming the variational Hodge
conjecture)
 of the same method to
the Chow motive of low degree $K3$ surfaces.

 Bloch's conjecture for
$0$-cycles on surfaces states the following:
\begin{conj} \label{blochconj}
(cf. \cite{blochbook}) Let $\Gamma\in CH^2(Y\times X)$, where $Y$ is
smooth projective and $X$ is a smooth projective surface. Assume
that $[\Gamma]^*:H^{2,0}(X)\rightarrow H^{2,0}(Y)$ vanishes. Then
$$\Gamma_*: CH_0(Y)_{alb}\rightarrow CH_0(X)_{alb}$$
vanishes.
\end{conj}
Here $[\Gamma]\in H^4(X\times X,\mathbb{Q})$ is the cohomology class
of $\Gamma$ and $$CH_0(Y)_{alb}:={\rm
Ker}\,(CH_0(Y)_{hom}\stackrel{alb_Y}{\rightarrow} Alb(Y)).$$

 Particular cases concern the case where $\Gamma$ is the
diagonal of  a surface $X$ with $q=p_g=0$. Then
$\Gamma_*=Id_{CH_0(X)_{alb}}$, so that the conjecture predicts
$CH_0(X)_{alb}=0$. This statement is known to hold for surfaces
which are not of general type by \cite{blochkas}, and for surfaces
of general type, it is known to hold by Kimura \cite{kimura} if the
surface $X$ is furthermore rationally dominated by a product of
curves (cf. \cite{BCGP} for many such examples). Furthermore, for
several other families of surfaces of general type with $q=p_g=0$,
it is known to hold either for the general member of the family (eg.
the Godeaux surfaces, cf. \cite{voisinscuola}), or for specific
members of the family (for example the Barlow surface
\cite{barlow2}).

A slightly more general situation concerns surfaces equipped with
the action of a finite group $G$. This has been considered in the
paper \cite{genhodgebloch}, where the following theorem concerning
group actions on complete intersection surfaces is proved: Let $X$
be a smooth projective variety with trivial Chow groups (i.e. the
cycle class map $CH_i(X)_\mathbb{Q}\rightarrow
H^{2n-2i}(X,\mathbb{Q}),\,n={\rm dim}\,X$ is injective for all $i$).
Let $G$ be a finite group acting on $X$ and let $E$ be a
$G$-equivariant rank $n-2$ vector bundle on $X$ which has ``enough''
$G$-invariant sections (for example, if the group action is trivial,
one asks that $E$ is very ample). Let $\pi\in \mathbb{Q}[G]$ be a
projector. Then $\pi$ gives a self-correspondence $\Gamma^\pi$ with
$\mathbb{Q}$-coefficients (which is a projector) of the
$G$-invariant surfaces $S=V(\sigma),\,\sigma\in H^0(X,E)^G$ (cf.
Section \ref{sec1}). We use the notation $$H^{2,0}(S)^\pi:={\rm
Im}\,([\Gamma^\pi]^*:H^{2,0}(S)\rightarrow H^{2,0}(S)),$$
$$CH_0(S)_{\mathbb{Q},hom}^\pi:={\rm
Im}\,([\Gamma^\pi]_*:CH_0(S)_{\mathbb{Q},hom}\rightarrow
CH_0(S)_{\mathbb{Q},hom})).$$
\begin{theo}\label{voisinhodgebloch} Assume that the smooth surfaces
$S=V(\sigma),\,\sigma\in H^0(E)^G$, satisfy $H^{2,0}(S)^\pi=0$. Then
we have $CH_0(S)_{\mathbb{Q},hom}^\pi=0$.
\end{theo}

 Note that the Bloch conjecture for finite group actions on surfaces
 which do not fit at all in the above geometric setting,
  namely finite order symplectic automorphisms of $K3$
surfaces, has been recently proved in \cite{huybrechts} and
\cite{voisininvok3}
 by completely different methods. For these
symplectic automorphisms, one considers the cycle
$\Delta_X-\frac{1}{|G|}\sum_{g\in G}{\rm Graph}\,g$, which acts as
the identity minus the projector onto the $G$-invariant part, and
one  proves that it acts  as $0$ on $CH_0(X)_{hom}$ (in fact on the
whole of $CH_0$) according to Conjecture \ref{blochconj}.

Theorem \ref{voisinhodgebloch} is  rather restrictive geometrically,
due to the fact that not only we consider $0$-sets of sections of a
vector bundle, but also  we impose this very ampleness assumption on
the vector bundle. Our first result in this paper is a relaxed
version of this theorem, which works in a much more general
geometric context and will be  applicable in particular to the case
of Catanese surfaces.

 Let  $\mathcal{S}\rightarrow B$ be a smooth projective morphism
with two dimensional connected fibers, where $B$ is
quasi-projective. Let $\Gamma\in
CH^2(\mathcal{S}\times_B\mathcal{S})_\mathbb{Q}$ be a relative
$0$-self correspondence. Let $\Gamma_t:=\Gamma_{\mid
\mathcal{S}_t\times\mathcal{S}_t}$ be the restricted cycle, with
cohomology class $[\Gamma_t]\in
H^4(\mathcal{S}_t\times\mathcal{S}_t,\mathbb{Q})$. We have the
actions
$$\Gamma_{t*}:CH_0(\mathcal{S}_t)_\mathbb{Q}\rightarrow CH_0(\mathcal{S}_t)_\mathbb{Q},\,\,
[\Gamma_t]^*:H^{i,0}(\mathcal{S}_t)\rightarrow
H^{i,0}(\mathcal{S}_t).
$$

\begin{theo} \label{theomain}
Assume the following:

 (1) The fibers $\mathcal{S}_t$ satisfy $h^{1,0}(\mathcal{S}_t)=0$
 and $[\Gamma_t]^*:H^{2,0}(\mathcal{S}_t)\rightarrow H^{2,0}(\mathcal{S}_t)$ is equal to zero.

(2) A  smooth projective (equivalently any smooth projective)
completion  $\overline{\mathcal{S}\times_B\mathcal{S}}$ of the
fibered self-product $\mathcal{S}\times_B\mathcal{S}$ is rationally
connected.

Then
$\Gamma_{t*}:CH_0(\mathcal{S}_t)_{hom}\rightarrow CH_0(\mathcal{S}_t)_{hom}$
is nilpotent for any $t\in B$.

\end{theo}
This statement is both weaker and stronger than Theorem
\ref{voisinhodgebloch} since on the one hand, the conclusion only
states the nilpotence of $\Gamma_{t*}$, and not its vanishing, while
on the other hand the geometric context is much more flexible and
the assumption on the total space of the family is much weaker.

In fact the nilpotence property is sufficient to imply the vanishing
in a number of situations which we describe below. The first
situation is the case where we consider a family of surfaces with
$h^{2,0}=h^{1,0}=0$. Then we get the following consequence (the
Bloch conjecture for surfaces with $q=p_g=0$ under assumption (2)
below):
\begin{coro} \label{corointro} Let $\mathcal{S}\rightarrow B$ be a smooth projective morphism
with two dimensional connected fibers, where $B$ is
quasi-projective. Assume the following:

 (1) The fibers $\mathcal{S}_t$ satisfy $H^{1,0}(\mathcal{S}_t)=H^{2,0}(\mathcal{S}_t)=0$.

 (2) A  projective
completion (or, equivalently, any projective
completion) $\overline{\mathcal{S}\times_B\mathcal{S}}$ of
the fibered self-product $\mathcal{S}\times_B\mathcal{S}$ is rationally connected.

Then $CH_0(\mathcal{S}_t)_{hom}=0$ for any $t\in B$.
\end{coro}

We refer to Section \ref{sec1}, Theorem \ref{variant} for a useful
variant involving group actions, which allows to consider many more
situations (cf. \cite{genhodgebloch} and Section \ref{sec2} for
examples).
\begin{rema}{\rm The proof will show as well that we can replace assumption (2) in
these statements by the following one:}

(2') A (equivalently any) smooth  projective completion
$\overline{\mathcal{S}\times_B\mathcal{S}}$ of the fibered
self-product $\mathcal{S}\times_B\mathcal{S}$ has trivial $CH_0$
group.

{\rm However, it seems more natural to put a geometric assumption on
the total space since this is in practice much easier to check.}
\end{rema}
In the second section of this paper, we will apply these results to
prove  Bloch's conjecture \ref{blochconj} for  Catanese surfaces
(cf. \cite{catanese}, \cite{supino}, \cite{BBK}). Catanese surfaces
can be constructed starting from a $5\times5$ symmetric matrix
$M(a),\,a\in\mathbb{P}^{11}$, of linear forms on $\mathbb{P}^3$
satisfying certain conditions (cf. (\ref{matrix})) making their
discriminant invariant under the Godeaux action (\ref{action}) of
$\mathbb{Z}/5\mathbb{Z}$ on $\mathbb{P}^3$. The general quintic
surface $V(a)$ defined by the determinant of $M(a)$ has $20$ nodes
corresponding to the points $x\in \mathbb{P}^3$ where the matrix
$M(a,x)$ has rank $3$, and it  admits  a double cover $S(a)$  which
is \'{e}tale away from the nodes, and  to which the
$\mathbb{Z}/5\mathbb{Z}$-action lifts. Then the Catanese surface
$\Sigma(a)$ is the quotient of $S(a)$ by this lifted action.

 Catanese surfaces
have a $4$-dimensional moduli space. For our purpose, the geometry
of the explicit $11$-dimensional parameter space described in
\cite{supino} is in fact  more important than the structure of the
moduli space.
\begin{theo}\label{theocamedelli} Let $S$ be a  Catanese surface. Then
$CH_0(S)=\mathbb{Z}$.
\end{theo}
The starting point of this work was a question asked by the authors
of \cite{BBK}: They needed to know that the Bloch conjecture holds
for a simply connected surface with $p_g=0$ (eg a Barlow surface
\cite{barlow1}) and furthermore, they needed it for a general
deformation of this surface. The Bloch conjecture was proved by
Barlow \cite{barlow2} for some Barlow surfaces admitting an extra
group action allowing to play on group theoretic arguments as in
\cite{inosemizukami}, but it was not known for the general Barlow
surface.

 Theorem
\ref{theocamedelli} implies as well the Bloch conjecture for the
Barlow surfaces since the Barlow surfaces can be constructed as
quotients of certain Catanese surfaces admitting an extra
involution, namely the determinantal equation defining $V(a)$ has to
be invariant under the action of the dihedral group of order $10$
(cf. \cite{supino}). The Catanese surfaces appearing in this
construction of Barlow surfaces have only
   a $2$-dimensional moduli space.
It is interesting to note that we get the Bloch conjecture for the
general Barlow surface via the Bloch conjecture for the general
Catanese surface, but that our strategy does not work directly for
the Barlow surface, which has a too small parameter space (cf.
\cite{catanese}, \cite{supino}).

The third section of this paper applies Theorem \ref{theomain} to
prove a conditional result on the Chow motive of $K3$ surfaces which
can be realized as $0$-sets of sections of a vector bundle on a
rationally connected variety (cf. \cite{mukai1}, \cite{mukai2}).
Recall that the Kuga-Satake construction (cf. \cite{kugasatake},
\cite{deligne}) associates to a polarized $K3$ surface $S$ an
abelian variety $K(S)$ with the property that the Hodge structure on
$H^2(S,\mathbb{Q})$ is a direct summand of the Hodge structure of
$H^2(K(S),\mathbb{Q})$. The Hodge conjecture predicts that the
corresponding degree $4$ Hodge class on $S\times K(S)$ is algebraic.
This is not known in general, but this is established for $K3$
surfaces with large Picard number (cf. \cite{morrison},
\cite{paranjape}). The next question concerns the Chow motive (as
opposed to the numerical motive) of these $K3$ surfaces. The
Kuga-Satake construction combined with the Bloch conjecture implies
that the Chow motive of a $K3$ surface is a direct summand of the
Chow motive of its Kuga-Satake variety. In this direction, we prove
the following Theorem \ref{theoK3}:
 Let $X$ be a rationally connected variety of dimension $n$ and let $E\rightarrow X$ be a
rank $n-2$ globally generated vector bundle satisfying  the
following properties:

(i)  The restriction map $H^0(X,E)\rightarrow H^0(z,E_{\mid z})$ is
surjective for general $z=\{x,\,y\}\subset X$.

(ii) The general section $\sigma$ vanishing at two general points
$x,\,y$ determines a smooth surface $V(\sigma)$.

 We consider the
case where
  the surfaces $S=V(\sigma)$ are  algebraic $K3$ surfaces.
  For example, this is the case if
 ${\rm det}\,E=-K_X$,  and the surfaces $S=V(\sigma)$ for general
$\sigma\in H^0(X,E)$ have irregularity $0$. Almost all general
algebraic $K3$ surfaces of genus $\leq 20$ have been described this
way  by Mukai (cf. \cite{mukai1}, \cite{mukai2}), where $X$ is a
homogeneous variety with Picard number $1$. Many more examples can
be constructed starting from an $X$ with Picard number $\geq 2$.
\begin{theo}\label{theoK3} Assume the variational Hodge conjecture in degree
$4$. Then the Chow motive of a $K3$ surface $S$ as above is a direct
summand of the Chow motive of an abelian variety.
\end{theo}
\begin{rema}{\rm The variational Hodge conjecture for degree $4$
Hodge classes is implied by the Lefschetz standard conjecture in
degrees $2$ and  $4$. It is used here only to conclude that the
Kuga-Satake correspondence is algebraic for any $S$ as above. Hence
we could replace the variational Hodge conjecture by the Lefschetz
standard conjecture or  by the assumption that the cohomological
motive of a general $K3$ surfaces $S$ in our family is a direct
summand of the cohomological motive of an abelian variety. The
contents of the theorem is  that we then have the same result  for
the Chow motive.}
\end{rema}
As a consequence of this result, we get the following (conditional) corollaries.
\begin{coro} \label{corok3}
With the same assumptions as in Theorem \ref{theoK3}, let $S$ be a
member of the family of $K3$ surfaces parameterized by
$\mathbb{P}(H^0(X,E))$, and let $\Gamma\in CH^2(S\times S)$ be a
correspondence such that $[\Gamma]^*:H^{2,0}(S)\rightarrow
H^{2,0}(S)$ is zero. Then $Z_*:CH_0(S)_{hom}\rightarrow
CH_0(S)_{hom}$ is nilpotent.
\end{coro}
\begin{rema}{\rm Note that there is a crucial difference between
Theorem \ref{theomain} and Corollary \ref{corok3}: In Corollary
\ref{corok3}, the cycle $\Gamma$ is not supposed to exist on the
general deformation $S_t$ of $S$. (Note also that the result in
Corollary \ref{corok3} is only conditional since we need the
Lefschetz standard conjecture, or at least to know that the
Kuga-Satake correspondence is algebraic for general $S_t$, while
Theorem \ref{theomain} is unconditional!)}
\end{rema}
\begin{coro}\label{corofinal} With the same assumptions as in Theorem \ref{theoK3}, the transcendental part
of the Chow motive of any member of the family of $K3$ surfaces
parameterized by $\mathbb{P}(H^0(X,E))$ is indecomposable, that is,
any submotive of it is either the whole motive or the $0$-motive.
\end{coro}
{\bf Thanks.}  I thank Christian B\"{o}hning, Hans-Christian Graf
von Bothmer, Ludmil Katzarkov and Pavel
 Sosna for asking me  the question whether  general Barlow surfaces satisfy
 the  Bloch conjecture,
 and for providing references on Barlow versus Catanese surfaces.

\section{Proof of Theorem \ref{theomain} and some consequences\label{sec1}}
This section is devoted to the proof of Theorem \ref{theomain} and
its consequences (Corollary \ref{corointro} or its more general form
Theorem \ref{variant}).

The proof will follow essentially the idea of  \cite{genhodgebloch}.
The main novelty in the proof lies in the use of Proposition
\ref{prop}. For completeness, we also outline the arguments of
\cite{genhodgebloch}, restricted to the surface case.

Consider the  codimension $2$-cycle
$$\Gamma\in CH^2(\mathcal{S}\times_B\mathcal{S})_\mathbb{Q}.$$

Assumption (1) tells us that the restricted cycle
$\Gamma_{t}:=\Gamma_{\mid \mathcal{S}_t\times\mathcal{S}_t}$ is
cohomologous to the sum of a cycle supported on a product of
(nonnecessarily irreducible curves) in $\mathcal{S}_t$ and of cycles
pulled-back from $CH_0(\mathcal{S}_t)$ via the two projections. We
deduce from this (cf. \cite[Prop. 2.7]{genhodgebloch}):
\begin{lemm} There exist a codimension $1$ closed algebraic subset
$\mathcal{C}\subset \mathcal{S}$,
 a codimension $2$ cycle $\mathcal{Z}$ on $\mathcal{S}\times_B\mathcal{S}$ with $\mathbb{Q}$-coefficients
  supported on
$\mathcal{C}\times_B\mathcal{C}$, and two codimension
 $2$ cycles $\mathcal{Z}_1,\,\mathcal{Z}_2$ with $\mathbb{Q}$-coefficients
 on
 $\mathcal{S}$, such that
the cycle
$$\Gamma-\mathcal{Z}-p_1^*\mathcal{Z}_1-p_2^*\mathcal{Z}_2$$ has
its restriction to each fiber $\mathcal{S}_t\times\mathcal{S}_t$
cohomologous to $0$, where
$p_1,\,p_2:\mathcal{S}\times_B\mathcal{S}\rightarrow\mathcal{S}$ are
the two projections.
\end{lemm}
This lemma is one of the key observations in \cite{genhodgebloch}.
The existence of the data above is rather clear  after a generically
finite base change $B'\rightarrow B$ because it is true fiberwise.
The key point is that, working with cycles with
$\mathbb{Q}$-coefficients, we can descend to $B$ and hence, do not
need in fact this base change which would ruin assumption (2).

  The next
step consists in passing from the fiberwise cohomological equality
$$[\Gamma-\mathcal{Z}-p_1^*\mathcal{Z}_1-p_2^*\mathcal{Z}_2]_{\mid \mathcal{S}_t\times\mathcal{S}_t}=0\,\,{\rm in}\,\,H^4(
\mathcal{S}_t\times\mathcal{S}_t,\mathbb{Q})$$
to the following global vanishing:
\begin{lemm}\label{secondlemme} (cf. \cite[Lemma 2.12]{genhodgebloch}) There exist codimension
$2$ algebraic cycles
$\mathcal{Z}'_1$, $\mathcal{Z}'_2$ with $\mathbb{Q}$-coefficients
on $\mathcal{S}$ such that
$$[\Gamma-\mathcal{Z}-p_1^*\mathcal{Z}'_1-p_2^*\mathcal{Z}'_2]=0\,\,{\rm in}\,\,H^4(
\mathcal{S}\times_B\mathcal{S},\mathbb{Q}).$$
\end{lemm}

The proof of this lemma consists in the study of the Leray spectral
sequence of the fibration
$p:\mathcal{S}\times_B\mathcal{S}\rightarrow B$. We know that the
class $[\Gamma-\mathcal{Z}-p_1^*\mathcal{Z}_1-p_2^*\mathcal{Z}_2]$
vanishes in the Leray quotient $H^0(B,R^4p_*\mathbb{Q})$ of
$H^4(\mathcal{S}\times_B\mathcal{S},\mathbb{Q})$. It follows that it
is of the form $p_1^*\alpha_1+p_2^*\alpha_2$, for some rational
cohomology classes $\alpha_1,\,\alpha_2$ on $\mathcal{S}$. One then
proves that $\alpha_i$ can be chosen to be algebraic on
$\mathcal{S}$.

The new part of the argument appears in the following proposition:
\begin{prop}\label{prop}  Under  assumption (2),  the following hold :

(i)  The cycle
$\Gamma-\mathcal{Z}-p_1^*\mathcal{Z}'_1-p_2^*\mathcal{Z}'_2$ is
algebraically equivalent to $0$ on $\mathcal{S}\times_B\mathcal{S}$.

(ii) The restriction to the fibers
$\mathcal{S}_t\times\mathcal{S}_t$ of the codimension $2$ cycle
$\mathcal{Z}'=\Gamma-\mathcal{Z}-p_1^*\mathcal{Z}'_1-p_2^*\mathcal{Z}'_2$
is a nilpotent element (with respect to the composition of
self-correspondences) of
$CH^2(\mathcal{S}_t\times\mathcal{S}_t)_\mathbb{Q}$.
\end{prop}
{\bf Proof.} We work now with a smooth projective completion
$\overline{\mathcal{S}\times_B\mathcal{S}}$. Let
$D:=\overline{\mathcal{S}\times_B\mathcal{S}}\setminus
{\mathcal{S}\times_B\mathcal{S}}$ be the divisor at infinity. Let
$\widetilde{D}\stackrel{j}{\rightarrow}\overline{\mathcal{S}\times_B\mathcal{S}}$
be a desingularization of $D$. The codimension $2$ cycle
$\mathcal{Z}'$ extends to a cycle $\overline{\mathcal{Z}'}$ over
$\overline{\mathcal{S}\times_B\mathcal{S}}$. We know from Lemma
\ref{secondlemme} that $$[\overline{\mathcal{Z}'}]=0\,\,{\rm in}\,\,
H^4(\overline{\mathcal{S}\times_B\mathcal{S}},\mathbb{Q})$$ and this
implies by \cite[Prop. 3]{voisintorino} that there is a degree $2$
Hodge class $\alpha$ on $\widetilde{D}$ such that
$$j_*\alpha=[\overline{\mathcal{Z}'}]\,\,{\rm in}\,\,
H^4(\overline{\mathcal{S}\times_B\mathcal{S}},\mathbb{Q}).$$
By the Lefschetz theorem on $(1,1)$-classes, $\alpha$ is the class of a codimension
$1$ cycle $\mathcal{Z}''$ of $\widetilde{D}$ and
we conclude that
$$[\overline{\mathcal{Z}'}-j_*\mathcal{Z}'']=0\,\,{\rm in}\,\,
H^4(\overline{\mathcal{S}\times_B\mathcal{S}},\mathbb{Q}).$$

We use now assumption (2) which says that the variety
$\overline{\mathcal{S}\times_B\mathcal{S}}$ is rationally connected.
It has then trivial $CH_0$, and so any codimension $2$ cycle
homologous to $0$ on $\overline{\mathcal{S}\times_B\mathcal{S}}$ is
algebraically equivalent to $0$ by the following result to to Bloch
and Srinivas:
\begin{theo}\cite{blochsrinivas} On a smooth projective variety $X$
with $CH_0(X)$ supported on a surface, homological equivalence and
algebraic equivalence coincide for codimension $2$ cycles.
\end{theo}
We thus conclude that $\overline{\mathcal{Z}'}-j_*\mathcal{Z}''$ is
algebraically equivalent to $0$ on
$\overline{\mathcal{S}\times_B\mathcal{S}}$, hence that
$\mathcal{Z}'=(\overline{\mathcal{Z}'}-j_*\mathcal{Z}'')_{\mid
\mathcal{S}\times_B\mathcal{S}}$ is algebraically equivalent to $0$
on ${\mathcal{S}\times_B\mathcal{S}}$.

(ii) This is a direct consequence of (i), using  the following
nilpotence result proved independently in \cite{voe} and
\cite{voisinsymmetric}:
\begin{theo} On any smooth projective variety,
self-correspondences algebraically equivalent to $0$ are nilpotent
for the composition of correspondences.
\end{theo}

 \cqfd
 {\bf Proof of Theorem \ref{theomain}.} Using the same notations as in the previous steps, we know by
Proposition  \ref{prop} that under assumptions (1) and (2), the
self-correspondence
\begin{eqnarray}
\label{douxnom}
\Gamma_{t}-\mathcal{Z}_t-p_1^*\mathcal{Z}'_{1,t}-p_2^*\mathcal{Z}'_{2,t}
\end{eqnarray}
on $\mathcal{S}_t$ with $\mathbb{Q}$-coefficients is nilpotent. In
particular, the morphism it induces at the level of Chow groups is
nilpotent. On the other hand, recall that $\mathcal{Z}_t$ is
supported on a product of curves on $\mathcal{S}_t$, hence acts
trivially on $CH_0(\mathcal{S}_t)_\mathbb{Q}$. Obviously, both
cycles $p_1^*\mathcal{Z}'_{1,t},\,p_2^*\mathcal{Z}'_{2,t}$
 act trivially on $CH_0(\mathcal{S}_t)_{\mathbb{Q},hom}$.
Hence the self-correspondence (\ref{douxnom}) acts as $\Gamma_{t*}$
on $CH_0(\mathcal{S}_t)_{\mathbb{Q},hom}$.
 \cqfd
  Let us now turn to our main application, namely Corollary \ref{corointro},
  or a more general form of it which involves a family of surfaces $S_t$ with an
action of a finite group $G$ and a
 projector $\pi\in \mathbb{Q}[G]$.
 Writing such a projector as $\pi=\sum_{g\in G}a_gg$, $a_g\in\mathbb{Q}$,
 such a projector provides a codimension $2$ cycle
 \begin{eqnarray}\label{gammapi}\Gamma^\pi_{t}=\sum_ga_g{\rm Graph}\,g\in CH^2(S_t\times S_t)_\mathbb{Q},
 \end{eqnarray}
  with actions
   $${\Gamma^\pi_{t}}^*=\sum_ga_gg^*$$
   on  the holomorphic forms of $S$, and
   $${\Gamma^\pi_{t}}_*=\sum_ga_gg_*$$
on $CH_0(S_t)_\mathbb{Q}$. Denote respectively
$CH_0(\mathcal{S}_t)_{\mathbb{Q},hom}^\pi$  the image of the
projector ${\Gamma^\pi_{t}}_*$ acting on
$CH_0(S_t)_{\mathbb{Q},hom}$ and $H^{2,0}(S_t)^\pi$ the image of the
projector ${\Gamma^\pi_{t}}^*$ acting on $H^{2,0}(S_t)$.
\begin{theo} \label{variant} Let $\mathcal{S}\rightarrow B$ be a smooth projective morphism
with two dimensional connected fibers, where $B$ is
quasi-projective. Let $G$ a finite group acting in a fiberwise way
on $\mathcal{S}$ and let $\pi\in \mathbb{Q}[G]$ be a projector.
Assume the following:

 (1) The fibers $\mathcal{S}_t$ satisfy $H^{1,0}(\mathcal{S}_t)=0$
 and $H^{2,0}(\mathcal{S}_t)^\pi=0$.

 (2) A  smooth projective
completion (or, equivalently, any smooth projective completion)
$\overline{\mathcal{S}\times_B\mathcal{S}}$ of the fibered
self-product $\mathcal{S}\times_B\mathcal{S}$ is rationally
connected.

Then $CH_0(\mathcal{S}_t)_{\mathbb{Q},hom}^\pi=0$ for any $t\in B$.
\end{theo}

{\bf Proof.} The group $G$ acts fiberwise on $\mathcal{S}\rightarrow
B$. Thus we have the universal cycle
$$\Gamma^\pi\in CH^2(\mathcal{S}\times_B\mathcal{S})_\mathbb{Q}$$
defined as $\sum_{g\in G}a_g{\rm Graph}\,g$, where the graph is
taken over $B$. Since by assumption the action of $[\Gamma^\pi]^*$
on $H^{2,0}(\mathcal{S}_t)$ is $0$, we can apply Theorem
\ref{theomain}, and conclude that $\Gamma^\pi_{t*}$ is nilpotent on
$CH_0(S_t)_{\mathbb{Q,}hom}$.
  On the other hand, $\Gamma^\pi_{t*}$ is  a projector onto
$CH_0(\mathcal{S}_t)_{\mathbb{Q},hom}^\pi$. The fact that it is
nilpotent implies thus that it is $0$, hence that
$CH_0(\mathcal{S}_t)_{\mathbb{Q},hom}^\pi=0$. \cqfd
\section{Catanese and Godeaux surfaces\label{sec2}}
Our main goal in this section is to check the main assumption of
Theorem \ref{variant}, namely the rational connectedness of the
fibered self-product of the universal  Catanese surface, in order to
prove Theorem \ref{theocamedelli}.

We follow \cite{catanese}, \cite{supino}. Catanese surfaces can be
described as follows: Consider the following symmetric $5\times 5$
matrix
\begin{eqnarray}\label{matrix}M_a=
\begin{pmatrix}
a_1x_1 & a_2x_2&a_3x_3&a_4x_4&0 \\
a_2x_2 & a_5x_3&a_6x_4&0&a_7x_1\\
a_3x_3&a_6x_4&0&a_8x_1&a_9x_2\\
a_4x_4&0&a_{8} x_1&a_{10}x_2&a_{11}x_3\\
0&a_7x_1&a_9x_2&a_{11}x_3&a_{12}x_4
\end{pmatrix}
\end{eqnarray}
depending on  $12$ parameters $a_1,\ldots, a_{12}$ and defining a
symmetric bilinear (or quadratic) form $q(a,x)$ on $\mathbb{C}^5$
depending on $x\in \mathbb{P}^3$.
 This is a homogeneous degree $1$ matrix in the variables
$x_1,\ldots, x_4$, and the vanishing of its determinant gives a
degree $5$ surface $V(a)$ in $\mathbb{P}^3$ with nodes at those
points $x=(x_1,\ldots,x_4)$  where the matrix has rank only $3$. We
will denote by $T$ the vector space generated by $a_1,\ldots,a_{12}$
and $B\subset \mathbb{P}(T)$ the open set of parameters $a$
satisfying this last condition.

The surface $V(a)$ is invariant under the Godeaux action of
$\mathbb{Z}/5\mathbb{Z}$ on $\mathbb{P}^3$, given by
$$g^*x_i=\zeta^ix_i,\,i=1,\ldots,\,4,$$
where $g$ is a generator of $\mathbb{Z}/5\mathbb{Z}$ and $\zeta$ is
a fifth primitive root of unity. This follows either from the
explicit development of the determinant (see \cite{supino}, where
one monomial is incorrectly written: $x_1x_2^2x_4^2$  should be
$x_1x_3^2x_4^2$) or from the following argument that we will need
later on: Consider the following linear action of
$\mathbb{Z}/5\mathbb{Z}$ on $\mathbb{C}^5$:
\begin{eqnarray}
\label{action} g(y_1,\ldots,y_5)=(\zeta
y_1,\zeta^2y_2,\ldots,\zeta^5y_5). \end{eqnarray}
 Then one checks
immediately that
\begin{eqnarray}\label{formula} q(a,x)(gy,gy')=\zeta q(a,gx)(y,y'),
\end{eqnarray}
 so that the discriminant of $q(a,x)$,
as a function of $x$, is invariant under the action of $g$.

The Catanese surface $\Sigma(a)$ is obtained as follows: there is a
natural  double cover $S(a)$ of $V(a)$, which is \'etale away from
the nodes, and parameterizes the rulings in the rank four quadric
$Q(a,x)$ defined by $q(a,x)$ for $x\in V(a)$. This double cover is
naturally equipped with a  lift of the
$\mathbb{Z}/5\mathbb{Z}$-action on $V(a)$, which is explicitly
described as follows:  Note first of all that the quadrics $Q(a,x)$
pass through the point $y_0=(0,0,1,0,0)$ of $\mathbb{P}^4$.
\begin{lemm} \label{petitlemme} For the general point $(a,x)\in
\mathbb{P}(T)\times\mathbb{P}^3$ such that $x\in V(a)$, the quadric
$Q(a,x)$ is not singular at the point $y_0$. \end{lemm}

{\bf Proof.} It suffices to exhibit one pair $(a,x)$ satisfying this
condition, and such that the surface $V(a)$ is well-defined, that is
the discriminant of $q(a,x)$, seen as a function of $x$, is not
identically $0$. Indeed, the family of surfaces $V(a)$ is flat over
the base near such a point, and the result for the generic pair
$(a,x)$ will then follow because  the considered property is open on
the total space of this family.

 We choose $(a,x)$ in such a way that the
first column of the matrix (\ref{matrix}) is $0$, so that $x\in
V(a)$. For example, we impose the conditions:
$$a_1=0,\,a_2=0,\,x_3=0,\,x_4=0.$$
Then the quadratic form $q(a,x)$ has for matrix $$
\begin{pmatrix}
0 & 0&0&0&0 \\
0 & 0&0&0&a_7x_1\\
0&0&0&a_8x_1&a_9x_2\\
0&0&a_{8} x_1&a_{10}x_2&0\\
0&a_7x_1&a_9x_2&0&0
\end{pmatrix}
$$
It is clear that $y_0$ is not a singular point of the corresponding
quadric if $a_8x_1\not=0$. On the other hand, for $a$ satisfying
$a_1=0,\,a_2=0$, and for a general point
$x=(x_1,\ldots,x_4)$, the matrix of $q(a,x)$ takes the form
$$
\begin{pmatrix}
0 & 0&a_3x_3&a_4x_4&0 \\
0 & a_5x_3&a_6x_4&0&a_7x_1\\
a_3x_3&a_6x_4&0&a_8x_1&a_9x_2\\
a_4x_4&0&a_{8} x_1&a_{10}x_2&a_{11}x_3\\
0&a_7x_1&a_9x_2&a_{11}x_3&a_{12}x_4.
\end{pmatrix}
$$
It is elementary to check that this matrix is generically of
maximal rank.
 \cqfd

It follows from this lemma that  for generic $a$ and generic  $x\in
V(a)$, the two rulings of $Q(a,x)$ correspond bijectively to the two
planes through $y_0$ contained in $Q(a,x)$. Hence for generic $a$,
the surface $S(a)$ is birationally equivalent to the surface $S'(a)$
parameterizing planes passing through $y_0$ and contained in
$Q(a,x)$ for some $x\in \mathbb{P}^3$ (which lies then necessarily
in $V(a)$). It follows from formula (\ref{formula}) that $S'(a)$,
which is a surface contained  in the Grassmannian $G_{0}$ of planes
in $\mathbb{P}^4$ passing through $y_0$, is invariant under the
$\mathbb{Z}/5\mathbb{Z}$-action on $G_0$ induced by (\ref{action}).
Hence there is a canonical $\mathbb{Z}/5\mathbb{Z}$-action on
$S'(a)$ which is compatible with the $\mathbb{Z}/5\mathbb{Z}$-action
on $V(a)$, and this immediately implies that the latter lifts to an
action on  the surface $S(a)$ for any $a\in B$.

 The Catanese surface is
defined as the quotient of $S(a)$ by $\mathbb{Z}/5\mathbb{Z}$. In
the following, we are going to apply Theorem \ref{variant}, and will
thus work with the universal family $\mathcal{S}\rightarrow B$ of
double covers, with its $\mathbb{Z}/5\mathbb{Z}$-action defined
above, where $B\subset\mathbb{P}^{11}$ is the Zariski open set
parameterizing smooth surfaces $S(a)$.

We prove now:
\begin{prop}\label{propoRC} The universal family $\mathcal{S}\rightarrow B$ has the
property that the fibered self-product
$\mathcal{S}\times_B\mathcal{S}$ has a rationally connected smooth
projective compactification.
\end{prop}
{\bf Proof.} By the description given above, the family
$\mathcal{S}\rightarrow B$ of surfaces $S(a),\,a\in B$, maps
birationally to an irreducible component of the following variety
\begin{eqnarray}\label{desc}
W=\{(a,x,[P])\in \mathbb{P}(T)\times\mathbb{P}^3\times
G_0,\,q(a,x)_{\mid P}=0\},
\end{eqnarray}
by the rational map which to a general point $(a,x), x\in V(a)$ and
a choice of ruling in the quadric $Q(a,x)$ associates the unique
plane $P$ passing through $y_0$, contained in $Q(a,x)$ and belonging
to the chosen ruling. It follows that
$\mathcal{S}\times_B\mathcal{S}$ maps birationally by the same map
(which we will call $\Psi$)  onto an irreducible component $W_2^0$
of the following variety
\begin{eqnarray}\label{descprod}
W_2:=\{(a,x,y,[P],[P'])\in
\mathbb{P}(T)\times\mathbb{P}^3\times\mathbb{P}^3\times G_0\times
G_0,\,q(a,x)_{\mid P}=0,\,q(a,y)_{\mid P'}=0\}.
\end{eqnarray}
Let $\mathcal{E}$ be the vector bundle of rank $5$ on $ G_0$ whose
fiber at a point $[P]$ parameterizing a plane  $P\subset
\mathbb{P}^4$ passing through $y_0$ is the space
$$H^0(P,\mathcal{O}_P(2)\otimes\mathcal{I}_{y_0}).$$
The family of
quadrics $q(a,x)$ on $\mathbb{P}^4$ provides a  $12$ dimensional
linear space $T$ of sections of the bundle
$$\mathcal{F}_0:=pr_1^*\mathcal{O}_{\mathbb{P}^3}(1)\otimes
pr_3^*\mathcal{E}\oplus pr_2^*\mathcal{O}_{\mathbb{P}^3}(1)\otimes
pr_{4}^*\mathcal{E}$$ on
$$Y_0:=\mathbb{P}^3\times\mathbb{P}^3\times
G_0\times G_0,$$ where as usual the $pr_i$'s denote the various
projections from $Y_0$ to its factors. For a point $a\in T$, the
corresponding section of $\mathcal{F}_0$ is equal to $(q(a,x)_{\mid
P},q(a,x')_{\mid P'})$ at the point $(x,x',[P],[P'])$ of $Y_0$.
 Formula (\ref{descprod}) tells
us that $W_2$ is  the zero set of the corresponding universal
section of the bundle
$${pr'_{1}}^*\mathcal{O}_{\mathbb{P}(T)}(1)\otimes
{pr'_{2}}^*\mathcal{F}_0$$ on $\mathbb{P}(T)\times Y_0$, where
the $pr'_i$ are now the two natural projections from
$\mathbb{P}(T)\times Y_0$ to its summands. Note that $W_2^0$
 has dimension
$15$, hence has the expected codimension $10$, since ${\rm
dim}\,\mathbb{P}(T)\times Y_0=25$.

There is now a subtlety in our situation: It is  not hard to see
that $T$ generates generically the bundle $\mathcal{F}_0$ on $Y_0$.
Hence there is a ``main'' component of $W_2$ which is also of
dimension $15$, and is generically fibered into $\mathbb{P}^1$'s
over $Y_0$. This component is not equal to $W_2^0$ for the following
reason: If one takes two general planes $P,\,P'$ through $y_0$, and
two general points $x,\,x'\in\mathbb{P}^3$, the conditions that
$q(a,x)$ vanishes identically on $P$ and $q(a,x')$ vanishes
identically on $P'$ implies that the third column of the matrix
$M(a)$ is identically $0$, hence that the point $y_0$ generates in
fact the kernel of both matrices $M(a,x)$ and $M(a,x')$.  On the
other hand, by construction
 of the map $\Psi$, generically along ${\rm Im}\,\Psi\subset W_2^0$, the point $y_0$ is a smooth
point of the quadrics $Q(a,x)$ and $Q(a,x')$.

The following lemma describes the component $W_2^0$.

\begin{lemm}\label{lemmeaussiplustard}
Let $\Phi: W_2^0\rightarrow Y_0=\mathbb{P}^3\times
\mathbb{P}^3\times G_0\times G_0$ be the restriction to $W_2^0$ of
the second projection $\mathbb{P}(T)\times Y_0\rightarrow Y_0$.

(i)  The image ${\rm Im}\,\Phi$ is a hypersurface $\mathcal{D}$ which admits a rationally connected
desingularization.

(ii) The generic rank of the evaluation map restricted to
$\mathcal{D}$: $$T\otimes \mathcal{O}_{\mathcal{D}}\rightarrow
{\mathcal{F}_0}_{\mid \mathcal{D}}$$ is $9$. (Note that the generic
rank of the evaluation map  $T\otimes \mathcal{O}_{Y_0}\rightarrow
{\mathcal{F}_0}$ is $10$.)

\end{lemm}
{\bf Proof.}  We already mentioned that the map $\Phi\circ \Psi$
cannot be dominating. Let us explain more precisely why, as this
will provide the equation for $\mathcal{D}$: Let $[P],\,[P']\in G_0$
and $x,\,x'\in\mathbb{P}^3$. The condition that
$(x,x',P,P')\in\Phi(W_2^0)={\rm Im}\,(\Phi\circ\Psi)$ is that for
some $a\in \mathbb{P}(V)$,
\begin{eqnarray}
\label{equations} q(a,x)_{\mid P}=0,\,q(a,x')_{\mid P'}=0
\end{eqnarray}
and furthermore that $y_0$ is a smooth point of both quadrics
$Q(a,x)$ and $Q(a,x')$.

Let $e,\,f\in P$ be vectors such that $y_0,e,f$ form a basis of $P$,
and choose similarly $e',f'\in P'$ in order to get a basis of $P'$.
Then among equations (\ref{equations}), we get
\begin{eqnarray}
\label{sousequations}
q(a,x)(y_0,e)=0,\,q(a,x)(y_0,f)=0,\,q(a,x')(y_0,x')(y_0,e')=0,\,q(a,x')(y_0,f')=0.
\end{eqnarray}
For fixed $P,\,P',\,x,\,x'$, these equations are linear forms in the
variables $a_3,\,a_6,\,a_8,\,a_9$ and it is not hard to see that
they are linearly independent for a generic choice of
$P,\,P',\,x,\,x'$, so that (\ref{sousequations}) implies that
$a_3=a_6=a_8=a_9=0$. But then, looking at the matrix $M(a)$ of
(\ref{matrix}) we see  that $y_0$ is in the kernel of $Q(a,x)$ and
$Q(a,x')$. As already mentioned, the latter does not happen
generically along ${\rm Im}\,\Psi$, and we deduce that ${\rm
Im}\,(\Phi\circ \Psi)$ is contained in the hypersurface
$\mathcal{D}$ where the four linear forms (\ref{sousequations}) in
the four variables $a_3,\,a_6,\,a_8,\,a_9$ are not independent.

We claim that $\mathcal{D}$ is rationally connected. To see this we
first prove that it is irreducible, which is done by restricting the
equation $f$ of  $\mathcal{D}$  (which is the determinant of the
$(4\times 4)$ matrix whose columns are the linear forms
(\ref{sousequations}) written in the basis $a_3,\,a_6,\,a_8,\,a_9$),
to a subvariety $Z$ of $Y_0$ of the form $Z=\mathbb{P}^3\times
\mathbb{P}^3\times C\times C'\subset Y_0$, where $ C,\,C'$ are
curves in $G_0$.

Consider the following $1$ dimensional families
$C\cong\mathbb{P}^1$, $C'\cong\mathbb{P}^1$ of planes passing
through $y_0$:
$$P_{t}=<e_1,\lambda e_2+\mu
e_4,e_3>,\,t=(\lambda,\mu)\in \mathbb{P}^1\cong C,$$
$$P'_{t'}=<e_5,\lambda'e_2+\mu'e_4,e_3>,\,t'=(\lambda',\mu')\in \mathbb{P}^1\cong C'.$$
 The equations (\ref{sousequations}) restricted to the parameters $t,t',x,x'$ give  the following four
combinations of $a_3,\,a_6,\,a_8,\,a_9$ depending on
$\lambda,\,\mu,\,x,x'$:

$$ a_3x_3,\,\,\lambda a_6x_4+\mu
a_8x_1,\,\,\,a_9x'_2,\,\,\lambda'a_6x'_4+\mu'a_8x'_1.$$ Taking the
determinant of this family gives
 $$f_{\mid
Z}=x_3x'_2(\lambda \mu'x_4x'_1-\lambda'\mu x'_4x_1).$$ The
hypersurface in $Z=\mathbb{P}^3\times \mathbb{P}^3\times C\times C'$
defined by $f_{\mid Z}$ has three irreducible components, which
belong respectively to the linear systems
$$|pr_1^*\mathcal{O}_{\mathbb{P}^3}(1)|,\,|pr_2^*\mathcal{O}_{\mathbb{P}^3}(1)|,\,|pr_{1}^*\mathcal{O}_{\mathbb{P}^3}(1)\otimes
pr_2^*\mathcal{O}_{\mathbb{P}^3}(1)\otimes
pr_3^*\mathcal{O}_{\mathbb{P}^1}(1)\otimes
pr_4^*\mathcal{O}_{\mathbb{P}^3}(1)|,
$$
where the $pr_i$'s are now the projections from
$Z=\mathbb{P}^3\times \mathbb{P}^3\times \mathbb{P}^1\times
\mathbb{P}^1$ to its factors.

 As the restriction map from ${\rm Pic}\,Y_0$ to ${\rm
Pic}\,Z$ is injective, we see that if $\mathcal{D}$ was reducible,
it would have an  irreducible component  in one of the linear
systems
$$|p_1^*\mathcal{O}_{\mathbb{P}^3}(1)|,\,|p_2^*\mathcal{O}_{\mathbb{P}^3}(1)|,\,|p_{1}^*\mathcal{O}_{\mathbb{P}^3}(1)\otimes
p_2^*\mathcal{O}_{\mathbb{P}^3}(1)\otimes
p_3^*\mathcal{O}_{\mathbb{P}^1}(1)\otimes
p_4^*\mathcal{O}_{\mathbb{P}^3}(1)|,
$$
where the $p_i$'s are now the projections from
$Y_0=\mathbb{P}^3\times \mathbb{P}^3\times G_0\times G_0$ to its
factors. In particular, its restriction to either any
$\mathbb{P}^3\times \{x'\}\times\{[P]\}\times \{[P']\}$ or any
$\{x\}\times\mathbb{P}^3\times \{[P]\}\times \{[P']\}$ would have an
irreducible component of degree $1$.

That this is not the case is proved by considering now the case
where
$$P=<e_1+e_2,\,e_4+e_5,\,e_3>, \,\,P'=<e_1+e_4,\,e_2+e_5,\,e_3>.$$
Then the four linear forms in (\ref{sousequations}) become
$$a_3x_3+a_6x_4,\,a_8x_1+a_9x_2,\,a_3x'_3+a_8x'_1,\,a_6x'_4+a_9x'_2.$$
It is immediate to check that for generic choice of $x'$, the
quadratic equation in $x$ given by this determinant has rank at
least $4$, and in particular is irreducible. Hence the restriction
of $f$ to $\mathbb{P}^3\times \{x'\}\times\{[P]\}\times \{[P']\}$
 has no irreducible component of degree $1$, and similarly for $\{x\}\times\mathbb{P}^3\times \{[P]\}\times
\{[P']\}$.

 Rational connectedness of $\mathcal{D}$ now follows
from the fact  that the  projection on the last three summands of
$Y_0$, restricted to $\mathcal{D}$,
$$p: \mathcal{D}\hookrightarrow Y_0\rightarrow \mathbb{P}^3\times G_0\times G_0,$$
 has for general fibers
smooth two-dimensional quadrics, as shows the computation we just
made. Hence $\mathcal{D}$ admits a rationally connected smooth
projective model, for example by \cite{GHS}.

We already proved that ${\rm Im}\,(\Phi\circ \Psi)\subset
\mathcal{D}$. We will now prove that ${\rm Im}\,(\Phi\circ \Psi)$
contains a Zariski open subset of $\mathcal{D}$ and  also statement
ii). This is done by the following argument which involves an
explicit computation at a point $a_0\in\mathbb{P}(T)$ where the
surface $V(a)$ becomes very singular but has enough smooth points to
conclude.

We start from the  following specific matrix :

\begin{eqnarray}\label{matrix0}
M_0=\begin{pmatrix}
x_1 & 0&x_3&0&0 \\
0 & x_3&x_4&0&0\\
x_3&x_4&0&x_1&0\\
0&0&x_1&x_2&0\\
0&0&0&0&x_4
\end{pmatrix}
\end{eqnarray}

We compute ${\rm det}\,M_0= -x_4x_3x_1^3-x_1x_2x_4^3-x_4x_2x_3^3 $.
The surface $V(a_0) $ defined by ${\rm det}\,M_0$ is rational and
smooth at the points
\begin{eqnarray}\label{x}x=(1,-\frac{1}{2},1,1),\,\,\,x'=(1,-\frac{2}{9},2,-1).
\end{eqnarray}
For each of the corresponding planes
$ P$ and $P'$  contained in $Q(a_0,x)$, $Q(a_0,x')$ respectively
and passing through $y_0$, we have two choices. We
choose the following:

\begin{eqnarray}\label{P}
P:=<y_0,\,e,\,f>,\,e=e_1+e_2-2e_4,\,f=e_2-e_4+\sqrt{-\frac{1}{2}}e_5,\\
\nonumber
P'=<y_0,\,e',\,f>,\,e'=4e_1-e_2-9e_4,\,f'=e_2+e_4+\frac{4}{3}e_5.
\end{eqnarray}
(We use here the standard basis $(e_1,\ldots,e_5)$ of
$\mathbb{C}^5$, so that $y_0=e_3$.) The vector $e-e_3$ (resp.
$e'-2e_3$) generates the kernel of $q(a_0,x)$ (resp. $q(a_0,x')$).
As   the quadrics $q(a_0,x)$ and $q(a_0,x')$ have rank $3$ and the
surface $V(a_0)$ is smooth at $x$ and $x'$, the points
 $(x,[P])$, $(x',[P'])$ determine smooth points of the surface
 $S(a_0)$.
 Even if the surface $S(a_0)$ is not smooth, so that $a_0\not\in B$,
 the universal family $\mathcal{S}\rightarrow B$ extends to a smooth
 non proper map $p:\mathcal{S}_e\rightarrow B_e\subset \mathbb{P}(T)$ near the points
 $(a_0,x,[P])$ and $(a_0,x',[P'])$.
Furthermore the morphism $\Psi$ is  well defined  at the points
$(a_0,x,P)$, $(a_0,x',P')$ since the point $y_0$ is not a singular
point of the quadric $Q(a_0,x)$ or $Q(a_0,x')$.

 We claim that the map
 $$\Phi\circ \Psi:\mathcal{S}_e\times_{B_e}\mathcal{S}_e\rightarrow Y_0=\mathbb{P}^3\times
 \mathbb{P}^3\times G_0\times G_0$$ has constant rank
 $13$ near the point
 $$((a_0,x,[P]), (a_0,x',[P']))\in \mathcal{S}_e\times_{B_e}\mathcal{S}_e.$$
 This implies that the image of $\Phi\circ \Psi$ is of dimension
 $13$, which is the dimension of
 $\mathcal{D}$, so that ${\rm Im}\,\Phi\circ \Psi$ has to contain a
 Zariski open set of $\mathcal{D}$ since $\mathcal{D}$ is
 irreducible. This implies that ${\rm Im}\,\Phi=\mathcal{D}$ as desired.

To prove the claim, we observe that it suffices to show the weaker
claim that the rank of \begin{eqnarray} \label{rank}(\Phi\circ
\Psi)_*:T_{\mathcal{S}_e\times_{B_e}\mathcal{S}_e}\rightarrow
T_{Y_0}
\end{eqnarray}
is equal to $13$ at the given point
 $((a_0,x,[P]), (a_0,x',[P']))$.
Indeed, as $\mathcal{S}_e\times_{B_e}\mathcal{S}_e$ is smooth
 at the point $((a_0,x,[P]), (a_0,x',[P']))$,
 the rank of the differential in (\ref{rank}) can only increase in a neighborhood of
 $((a_0,x,[P]), (a_0,x',[P']))$ and on the other hand, we know that it
 is not of maximal rank at any point of $\mathcal{S}'\times_{B'}\mathcal{S}'$,
 since ${\rm Im}\,(\Phi\circ \Psi)\subset \mathcal{D}$. Hence it must stay constant
 near $((a_0,x,[P]), (a_0,x',[P']))$.

 We now compute the rank of $(\Phi\circ \Psi)_*$  at the point
 $((a_0,x,[P]), (a_0,x',[P']))$.
  Note that the birational map
 $\Psi: \mathcal{S}_e\times_{B_e}\mathcal{S}_e \rightarrow W_2$ is  a local isomorphism
 near the point $((a_0,x,[P]), (a_0,x',[P']))$, because
 it is well defined at this point and its inverse too. This argument
 proves that not only $W_2^0$ but also $W_2$ is smooth of dimension
 $15$ near the point $(a_0,x,x',[P],[P'])$.
 Hence it suffices to compute the rank
 of the differential
 \begin{eqnarray}\label{diff}\Phi_*:T_{W_2,(a_0,x,x',[P],[P'])}\rightarrow T_{Y_0,(x,x',[P],[P'])}
 \end{eqnarray} of
 the map
 $\Phi: W_2\rightarrow Y_0$
 at the point $(a_0,x,x',[P],[P'])$.

Recalling that $W_2$ is defined as the zero
set of the universal section
$\sigma_{univ}$ of the bundle
$${pr'_{1}}^*\mathcal{O}_{\mathbb{P}(T)}(1)\otimes
{pr'_{2}}^*\mathcal{F}_0$$ on $\mathbb{P}(T)\times Y_0$,
the tangent space
of $W_2$ at the point
$(a_0,x,x',[P],[P'])$ is equal to
$${\rm Ker}\,(d\sigma_{univ}:T_{\mathbb{P}(T),a_0}\times T_{Y_0,(x,x',[P],[P'])}\rightarrow \mathcal{F}_{0,(x,x',[P],[P'])}).$$
Clearly the differential $d\sigma_{univ}$ restricted to
$T_{\mathbb{P}(T),a_0}$ is induced by the evaluation map
\begin{eqnarray}
\label{ev}{\rm ev}_{(x,x',[P],[P'])}:T\rightarrow
\mathcal{F}_{0,(x,x',[P],[P'])}
\end{eqnarray}
at the point $(x,x',[P],[P'])$ of $Y_0$. On the other hand, the
differential $d\sigma_{univ}$ restricted to the tangent space
$T_{Y_0,(x,x',[P],[P'])}$ is surjective because the variety
$S(a_0)\times S(a_0)$ is smooth of codimension $4$ and isomorphic
via $\Phi$ to $V(\sigma_{univ}(a_0))$ near $(x,x',[P],[P'])$. It
follows from this that the corank of the second projection
$$\Phi_*:{\rm Ker}\,d\sigma_{univ}=T_{W_2,(a_0,x,x',[P],[P'])}\subset T_{\mathbb{P}(T),a_0}\oplus T_{Y_0,(x,x',[P],[P'])}\rightarrow
T_{Y_0,(x,x',[P],[P'])},$$
  is equal to the corank of the map
  $d\sigma_{univ}:T_{\mathbb{P}(T),a_0}\rightarrow
  \mathcal{F}_{0,(x,x',[P],[P'])})$, that is to the corank of the
evaluation map ${\rm ev}_{(x,x',[P],[P'])}$ of (\ref{ev}).

In particular,  the rank of $\Phi_*$ in (\ref{diff})  is equal to
$13$ if and only if the rank of the evaluation map ${\rm
ev}_{x,x',[P],[P']}$  is equal to $9$, which is our statement (ii).

In conclusion, we proved that (i) is implied by (ii) and that (ii)
itself implied by (ii) at the given point $(x,x',[P],[P'])$ of
(\ref{x}), (\ref{P}).
 It just remains to prove that the rank of ${\rm ev}_{(x,x',[P],[P'])}$ is equal to $9$ at this point, which
is done by the explicit computation of the rank of the family of
linear forms in the $a_i$ given by
$$q(a,x)(y_0,e),\,\,\,q(a,x)(y_0,f),\,\,\,q(a,x)(e,e),\,\,\,q(a,x)(f,f),\,\,\,q(a,x)(e,f),$$
$$q(a,x')(y_0,e'),\,\,\,q(a,x')(y_0,f'),\,\,\,q(a,x')(e',e'),\,\,\,q(a,x')(f',f'),\,\,\,q(a,x')(e',f').$$
These forms are the following: $$ a_3+a_6-2a_8,\,\,\,\,\,\,
a_6-a_8-\frac{1}{2}\sqrt{-\frac{1}{2}}a_9,\,\,\,\,\,\,
a_1+a_5-2a_{10}-a_2-4a_4,$$
 $$
a_5-\frac{a_{10}}{2}-\frac{a_{12}}{2}+2\sqrt{-\frac{1}{2}}a_7
-2\sqrt{-\frac{1}{2}}a_{11},\,\,\,\,\,\,-\frac{a_{2}}{2}-a_4+a_5+\sqrt{-\frac{1}{2}}a_7-a_{10}-2\sqrt{-\frac{1}{2}}a_{11}
$$
$$ 8a_3+a_6-9a_8,\,\,\,\,\,\,-a_6+a_8-\frac{8}{27}a_9,\,\,\,\,\,\,16 a_1+2a_5-18
a_{10}+\frac{16}{9}a_2+72 a_4,$$
 $$
2a_5-\frac{2}{9}a_{10}-\frac{16}{9}a_{12}+\frac{8}{3}a_7+\frac{16}{3}a_{11},\,\,\,-\frac{8}{9}a_2-4a_4-2a_5-\frac{4}{3}a_7+2a_{10}-24
a_{11}.
$$

It is immediate to check that the rank of this family is  $9$.

 \cqfd

 It follows from Lemma \ref{lemmeaussiplustard}
 that  $W_2^0$ (or rather a smooth projective birational model
  of $W_2^0$)  is rationally connected, since
  a Zariski open set of $W_2^0$ is a $\mathbb{P}^2$-bundle over a Zariski open set
  $\mathcal{D}^0$ of $\mathcal{D}$ which is smooth
   and  admits a rationally connected completion. Hence
 $\mathcal{S}\times_B\mathcal{S}$ is also rationally connected,
as it is birationally equivalent  to $W_2^0$.
 \cqfd
 We get by application of Theorem
 \ref{variant} the following statement (cf. Theorem \ref{theocamedelli}):
\begin{coro}\label{corocamp} The Catanese surface $\Sigma(a)$ has $CH_0$ equal to
$\mathbb{Z}$.
\end{coro}
{\bf Proof.} Indeed, the surface $\Sigma(a)$ has $h^{1,0}=h^{2,0}=0$
(cf. \cite{catanese})  which means equivalently that the surface
$S(a)$ introduced above has $h^{2,0}_{inv}=0$ where ``inv'' means
invariant under the $\mathbb{Z}/5\mathbb{Z}$-action, which by the
above description is defined  on $\mathcal{S}$. The result  is then
a consequence of Theorem \ref{variant} applied to the projector
$\pi_{inv}$ onto the $\mathbb{Z}/5\mathbb{Z}$-invariant part, which
shows that
$$CH_0(\Sigma(a))_{\mathbb{Q},hom}=CH_0(S(a))^{\pi_{inv}}_{\mathbb{Q},hom}=0$$
and from Roitman's theorem \cite{roitman} which says that
$CH_0(\Sigma(a))$ has no torsion.
 \cqfd
 \begin{coro} The Barlow surface $\Sigma'(b)$  has
 $CH_0$ equal to $\mathbb{Z}$.
 \end{coro}
 {\bf Proof.} Indeed, the Barlow surface $\Sigma'(b)$ is a quotient
 of the Catanese surface $\Sigma(b)$ by an involution, hence
 $CH_0(\Sigma'(b))\hookrightarrow CH_0(\Sigma(b))$ since
 $CH_0(\Sigma'(b))$ has no torsion by \cite{roitman}.
 \cqfd
\section{On the Chow motive of complete intersections $K3$ surfaces}

Let $S$ be a $K3$ surface. The Hodge structure on $H^2(S,\mathbb{Z})$ is a weight $2$ polarized Hodge structure with $h^{2,0}=1$.
In \cite{kugasatake}, Kuga and Satake construct an abelian variety
$K(S)$ associated
to this Hodge structure. Its main property is the fact that
there is a natural injective morphism of weight $2$ Hodge structures :
$$H^2(S,\mathbb{Z})\rightarrow H^2(K(S),\mathbb{Z}).$$

Such a morphism of Hodge structures in turn provides, using K\"{u}nneth decomposition
 and Poincar\'e duality a Hodge class
$$\alpha_S\in {\rm Hdg}^4(S\times K(S))$$
where ${\rm Hdg}^4(X)$ denotes the space of rational Hodge classes on $X$.

This class is not known in general  to satisfy the Hodge conjecture, that is, to be the class of an algebraic cycle.
This is known to hold for Kummer surfaces (see \cite{morrison}), and for some
$K3$ surfaces with Picard number $16$ (see \cite{paranjape}).
The deformation theory of $K3$ surfaces, and more particularly the fact
 that any  projective $K3$ surface  deforms to a Kummer surface,
 combined with the global invariant cycle theorem
 of Deligne \cite{delignetheorie}  imply the following result, (cf.
\cite{deligne}, \cite{andre2}):

  \begin{theo} \label{theogenstuf} Let $S$ be a projective $K3$ surface. There exist
 a connected quasi-projective variety $B$, a family of projective  $K3$ surfaces
 $\mathcal{S}\rightarrow B$,  a family of abelian varieties
 $\mathcal{K}\rightarrow B$, (where all varieties are quasi-projective and all morphisms
 are smooth projective), and a Hodge class
 $$\eta\in {\rm Hdg}^4(\mathcal{S}\times_B\mathcal{K}),$$
 such that :

 \begin{enumerate}
 \item \label{1}for some point $t_0\in B$, $\mathcal{S}_{t_0}\cong S$;
 \item \label{2}
 For any point $t\in B$, $\mathcal{K}_{t}\cong K(S_t)$ and
 $\eta_t=\alpha_{S_t}$.
 \item For some point $t_1\in B$, the class $\eta_{t_1}$ is algebraic.
 \end{enumerate}
 \end{theo}
 Here by ``Hodge class $\eta$'', we mean that the class $\eta$ comes from a
 Hodge class on some smooth projective compactification
 of $\mathcal{S}\times_B\mathcal{K}$.

 The existence of the family is an algebraicity statement for the Kuga-Satake construction
 (see \cite{deligne}), which can then be done in family
 for $K3$ surfaces with given polarization,
 while the last item follows from the fact that the locally complete such families always contain Kummer fibers
 $S'$
 for which the class $\alpha_{S'}$ is known to be algebraic.
 \begin{coro} \cite{andre2} The Hodge class $\alpha_S$ is ``motivated''.
 \end{coro}
This means  (cf. \cite{andre2}) that  this Hodge class can be
constructed via algebraic correspondences
 from Hodge classes on auxiliary varieties, which  are either algebraic or
 obtained by inverting Lefschetz operators. In particular the class $\alpha_S$ is algebraic if the standard Lefschetz conjecture holds.
 \begin{coro} Assume the variational form of the Hodge conjecture holds, or assume the Lefschetz standard conjecture in degrees $2$ and $4$. Then the class
 $\alpha_S$ is algebraic for any projective $K3$ surface $S$.

 \end{coro}
 Indeed, the variational Hodge conjecture states that in the situation
 of Theorem \ref{theogenstuf}, if a Hodge class $\eta$ on the total space
 has an algebraic  restriction on one fiber, then its restriction to {\it any} fiber is algebraic. In our situation, it will be implied by the Lefschetz conjecture
 for degree $2$ and degree $4$ cohomology on a smooth projective compactification of
  $\mathcal{S}\times_B\mathcal{K}$.

From now on, we assume that the Kuga-Satake correspondence
$\alpha_S$ is algebraic for a general projective $K3$ surface $S$ of
genus $g$ (which means by definition that $S$ comes equipped with an
ample line bundle  of self-intersection $2g-2$).  We view the class
$\alpha_S$ as an injective  morphism of Hodge structures
\begin{eqnarray}\label{morphism}\alpha_S:
H^2(S,\mathbb{Q})\rightarrow H^2(K(S),\mathbb{Q}) \end{eqnarray}
 and
use a polarization $h$ of $K(S)$ to construct an inverse of
$\alpha_S$ by the following lemma (which can be proved as well by
the explicit description $K(S)$ as a complex torus and its
polarization, cf. \cite{kugasatake}). In the following,
$H^2(S,\mathbb{Q})_{tr}$ denotes the transcendental cohomology of
$S$, which is defined as the orthogonal complement of the
N\'eron-Severi group of $S$.
\begin{lemm}\label{pol} Let $h=c_1(H)\in  H^2(K(S),\mathbb{Q})$ be the class of an ample line bundle, where $S$ is a very general $K3$ surface of genus $g$.
Then there exists a nonzero rational number
$\lambda_g$ such that the endomorphism
$$^t\alpha_S\circ (h^{N-2}\cup)\circ \alpha_S:H^2(S,\mathbb{Q})\rightarrow H^2(S,\mathbb{Q})$$
restricts  to $\lambda_g Id$ on $H^2(S,\mathbb{Q})_{tr}$, where $N={\rm dim}\,K(S)$ and
$$^t\alpha_S: H^{2N-2}(K(S),\mathbb{Q})\rightarrow H^2(S,\mathbb{Q})$$
is the transpose of the map $\alpha_S$ of (\ref{morphism}) with
respect to Poincar\'{e} duality.
\end{lemm}
{\bf Proof.} The composite $^t\alpha_S\circ h^{N-2}\circ \alpha_S$ is an endomorphism
of the Hodge structure on $H^2(S,\mathbb{Q})$. As $S$ is very general, this Hodge structure
has only a one dimensional
$\mathbb{Q}$-vector space of algebraic classes, generated by the polarization,
and its orthogonal is a simple Hodge structure with
only the homotheties as endomorphisms.
We conclude that $^t\alpha_S\circ h^{N-2}\circ \alpha_S$ preserves
$H^2(S,\mathbb{Q})_{tr}$ and acts on it as an homothety  with rational coefficient
(which is thus independent of $S$).
It just remains to show that it does not act as zero on
$H^2(S,\mathbb{Q})_{tr}$. This follows from the second Hodge-Riemmann
bilinear relations which say that for
$\omega\in H^{2,0}(S),\,\omega\not=0$, we have
$$\langle\omega, \,{^t\alpha_S}(h^{N-2}\cup\alpha_S(\overline{\omega}))\rangle_S= \langle\alpha_S(\omega),h^{N-2}\cup\alpha_S(\overline{\omega})\rangle_{K(S)}$$
which is $>0$ because $\alpha_S(\omega)\not=0$ in $H^{2,0}(K(S))$.
Hence $^t\alpha_S(h^{N-2}\cup\alpha_S(\overline{\omega}))\not=0$.

\cqfd
We now start from a rationally connected variety $X$ of dimension $n$, with a vector bundle
$E$ of rank $n-2$ on $X$, such that  $-K_X={\rm det}\,E$ and the following properties hold:

(*) {\it  For general $x,\,y\in X$, and for general $\sigma\in
H^0(X,E\otimes\mathcal{I}_{x}\otimes \mathcal{I}_y)$, the zero locus
$V(\sigma)$ is a smooth connected surface with $0$ irregularity
(hence a smooth $K3$ surface).}

Let $L$ be an ample line bundle on $X$, inducing a polarization
of genus $g$
on the $K3$ surface $S_\sigma:=V(\sigma)$ for
$\sigma\in B\subset \mathbb{P}(H^0(X,E))$.

We now prove:
\begin{theo} \label{theocond}
Assume the Kuga-Satake correspondence $\alpha_S$ is algebraic, for
the general $K3$ surface $S$ with such a polarization. Then, for any
$\sigma\in B$, the Chow motive of $S_\sigma$ is a direct summand of
the motive of an abelian variety.
\end{theo}
Let us  explain the precise meaning of this statement. The
algebraicity of the Kuga-Satake correspondence combined with Lemma
\ref{pol} implies that there is a codimension $2$ algebraic cycle
$$Z_S\in CH^2(S\times K(S))_\mathbb{Q}$$
with the property that
 the cycle $\Gamma_S$ defined by
$$ \Gamma_S={^tZ_S}\circ h^{N-2}\circ Z_S\in CH^2(S\times S)_\mathbb{Q}$$
has the property that its cohomology class
$[\Gamma]\in H^4(S\times S,\mathbb{Q})$ induces
a nonzero homothety
\begin{eqnarray}\label{coheqretour}[\Gamma]_*=\lambda Id:H^2(S,\mathbb{Q})_{tr}\rightarrow H^2(S,\mathbb{Q})_{tr},
\end{eqnarray}
which can be equivalently formulated as follows: Let us introduce
the cycle $\Delta_{S,tr}$, which in the case of a $K3$ surface is
canonically defined, by the formula
\begin{eqnarray}
\label{gogo}\Delta_{S,tr}=\Delta_S-o_S\times S-S\times
o_S-\sum_{ij}\alpha_{ij} C_i\times C_j, \end{eqnarray}
 where
$\Delta_S$ is the diagonal of $S$, $o_S$ is the canonical $0$-cycle
of degree $1$ on $S$ introduced in \cite{beauvillevoisin}, the $C_i$
form a basis of $({\rm Pic}\,S)\otimes\mathbb{Q}={\rm
NS}(S)_\mathbb{Q}$ and the $\alpha_{ij}$ are the coefficients of the
inverse of the matrix of the intersection form of $S$ restricted to
$NS(S)$. This corrected diagonal cycle is a projector and it has the
property that its action on cohomology is the orthogonal projector
$H^*(S,\mathbb{Q})\rightarrow H^2(S,\mathbb{Q})_{tr}$. Formula
(\ref{coheqretour}) says that we have the cohomological equality
\begin{eqnarray}\label{coheqretouravec deltatr}[\Gamma\circ \Delta_{S,tr}]=\lambda
[\Delta_{S,tr}]\,\,{\rm in}\,\,H^4(S\times S,\mathbb{Q}).
\end{eqnarray}
A more precise form of
Theorem \ref{theocond} says that we can get in fact such an equality at the level of
Chow groups:
\begin{theo} \label{autreforme}
Assume the Kuga-Satake correspondence $\alpha_S$ is algebraic, for a
general $K3$ surface with such a polarization. Then, for any
$\sigma\in B$, there is an abelian variety $A_\sigma$, and cycles
$Z\in CH^2(S_\sigma\times A_\sigma)_\mathbb{Q}$, $Z'\in
CH^{N'}(A_\sigma\times S_\sigma)_\mathbb{Q}$, $N'={\rm
dim}\,A_\sigma$, with the property that
\begin{eqnarray}\label{formuleaprouver} Z'\circ Z\circ \Delta_{S,tr}=\lambda
\Delta_{S,tr} \end{eqnarray} for a nonzero rational number
$\lambda$.
\end{theo}

The proof will use two preparatory lemmas.
 Let $B\subset \mathbb{P}(H^0(X,E))$ be the open set
parameterizing smooth surfaces
$S_\sigma=V(\sigma)\subset X$. Let
$\pi:\mathcal{S}\rightarrow B$ be the universal family, that is:
\begin{eqnarray}
\label{Sfin}\mathcal{S}=\{(\sigma,x)\in B\times X,\,\sigma(x)=0\}.
\end{eqnarray}
The first observation is the following:
\begin{lemm} \label{le1fin} Under assumption (*), the fibered self-product $\mathcal{S}\times_B\mathcal{S}$
is rationally connected (or rather, admits a smooth projective rationally connected
completion).
\end{lemm}
{\bf Proof.} From (\ref{Sfin}), we deduce
that
\begin{eqnarray}
\label{SxSfin}\mathcal{S}\times_B\mathcal{S}=\{(\sigma,x,y)\in B\times X\times X,\,\sigma(x)=\sigma(y)=0\}.
\end{eqnarray}
Hence $\mathcal{S}\times_B\mathcal{S}$ is Zariski open in
the following variety:
$$W:=\{(\sigma,x,y)\in \mathbb{P}(H^0(X,E))\times X\times X,\,\sigma(x)=\sigma(y)=0\}.$$
In particular, as it is irreducible, it is Zariski open in one irreducible component
$W^0$ of $W$.

Consider the projection on the two last factors:
$$(p_2,p_3): W\rightarrow X\times X.$$
Its fibers are projective spaces, so that there is only one ``main''
irreducible component $W^1$  of $W$ dominating $X\times X$ and it
admits a smooth
 rationally connected completion since $X\times X$ is rationally connected.

Assumption (*) now tells us that at a general point of $W^1$, the
first projection $p_1:W\rightarrow B$ is smooth of relative
dimension $4$. It follows that $W$ is smooth at this point which
belongs to both components $W^0$ and $W^1$. Thus $W^0=W^1$ and
$\mathcal{S}\times_B\mathcal{S}\cong_{birat} W^0$ admits a smooth
 rationally connected completion.
\cqfd
The next step is the following lemma:
\begin{lemm} \label{le2fin} Assume  the Kuga-Satake correspondence
$\alpha_S$ is algebraic for the general  polarized $K3$ surface of
genus $g$. Then there exist a rational number $\lambda\not=0$, a
family $\mathcal{A}\rightarrow B$ of polarized $N'$-dimensional
abelian varieties, with relative polarization $\mathcal{L}$ and a
codimension $2$ cycle
$$\mathcal{Z}\in CH^2(\mathcal{S}\times_B \mathcal{A})_\mathbb{Q}$$
such that for very general $t\in B$, the cycle
$$\Gamma_t:=   {^t\mathcal{Z}_t}\circ c_1(\mathcal{L}_t)^{N'-2}\circ \mathcal{Z}_t\circ \Delta_{tr,t}$$
satisfies:
\begin{eqnarray}\label{eqcohfin}[\Gamma_t]=\lambda[\Delta_{tr,t}]\in H^4(S_t\times S_t,\mathbb{Q}).
\end{eqnarray}

\end{lemm}
In this formula, the term $c_1(\mathcal{L}_t)^{N'-2}$ is defined as
the self-correspondence of $\mathcal{A}_t$ which consists of the
cycle $c_1(\mathcal{L}_t)^{N'-2}$ supported on the diagonal of
$\mathcal{A}_t$. We also recall  that the codimension $2$-cycle
$\Delta_{tr,t}\in CH^2(S_t\times S_t)_\mathbb{Q}$ is the projector
onto the transcendental part of the motive of $S_t$. The reason why
the result is stated only for the very general point $t$ is the fact
that due to the possible jump of the Picard group of $S_t$, the
generic cycle $\Delta_{tr,\overline{\eta}}$ does not specialize to
the cycle $\Delta_{tr,t}$ at any closed point $t\in B$, but only at
the very general one. In fact, the statement is true at any point,
but the cycle $\Delta_{tr,t}$ has to be modified when the Picard
group jumps.

\vspace{0.5cm}

{\bf Proof of Lemma \ref{le2fin}.} By our assumption, using the
countability of relative Hilbert schemes and the existence of
universal objects parameterized by them, there exist a generically
finite cover $r:B'\rightarrow B$, a universal family of polarized
abelian varieties
$$\mathcal{K}\rightarrow B',\,\mathcal{L}_K\in {\rm Pic}\,\mathcal{K}$$
and a codimension $2$ cycle $\mathcal{Z}'\in
CH^2(\mathcal{S}'\times_{B'}\mathcal{K})_\mathbb{Q}$, where
$\mathcal{S}':=\mathcal{S}\times_{B}B'$, with the property that
\begin{eqnarray}\label{autrebonne}
[\mathcal{Z}_t]=\alpha_{\mathcal{S}_t}\,\,{\rm
in}\,\,H^4(\mathcal{S}_t\times K(\mathcal{S}_t),\mathbb{Q})
\end{eqnarray}
 for any $t\in B'$. Furthermore,
by Lemma \ref{pol}, we know that there exists a nonzero rational
number $\lambda_g$ such that for any $t\in B'$, we have
\begin{eqnarray}
\label{labonne}^t\alpha_{\mathcal{S}_t}\circ (h^{N-2}\cup)\circ
\alpha_{\mathcal{S}_t}=\lambda_g
Id:H^2(\mathcal{S}_t,\mathbb{Q})_{tr}\rightarrow
H^2(\mathcal{S}_t,\mathbb{Q})_{tr}, \end{eqnarray}
 where as before $N$ is
the dimension of $K(\mathcal{S}_t)$, and $h_t=c_1(\mathcal{L}_{K\mid
\mathcal{K}_t})$. We now construct the following family of abelian
varieties on $B$ (or a Zariski open set of it)
$$\mathcal{A}_t=\prod_{t'\in r^{-1}(t)}\mathcal{K}_{t'},$$
with  polarization given by
\begin{eqnarray}\mathcal{L}_t=\sum_{t'\in
r^{-1}(t)}pr_{t'}^*(\mathcal{L}_{K,t'}), \label{polforfin}
\end{eqnarray} where $pr_{t'}$ is the obvious projection from
$\mathcal{A}_t=\prod_{t'\in r^{-1}(t)}K_{t'}$ to its factor
$\mathcal{K}_{t'}$, and the following cycle $\mathcal{Z}\in
CH^2(\mathcal{S}\times _B\mathcal{K})_\mathbb{Q}$, with fiber at
$t\in B$ given by
\begin{eqnarray}\label{cycleforfin}\mathcal{Z}_t=\sum_{t'\in
r^{-1}(t)}(Id_{\mathcal{S}_t},pr_{t'})^*\mathcal{Z}'_{t'}.
\end{eqnarray} In the
last formula, we use of course the identification
$$\mathcal{S}'_{t'}=\mathcal{S}_t,\,r(t')=t.$$
It just remains to prove formula (\ref{eqcohfin}). Combining
(\ref{polforfin}) and (\ref{cycleforfin}), we get, using again the
notation $h_t=c_1(\mathcal{L}_t)\in H^2(A_t,\mathbb{Q})$,
$h_{t'}=c_1(\mathcal{L}_{K,t'})\in H^2(K_{t'},\mathbb{Q})$:

$$[\Gamma_t]^*=$$
$$(\sum_{t'\in
r^{-1}(t)}[(pr_{t'},Id_{\mathcal{S}_t})^*(^t\mathcal{Z}'_{t'})]^*)\circ
(\sum_{t''\in r^{-1}(t)}pr_{t''}^*(h_{t''}))^{N'-2}\cup\circ
(\sum_{t'''\in
r^{-1}(t)}(Id_{\mathcal{S}_t},pr_{t'''})^*[\mathcal{Z}'_{t'''}]^* )
\circ \pi_{tr,t}:
$$
$$H^*(\mathcal{S}_t,\mathbb{Q})\rightarrow
H^*(\mathcal{S}_t,\mathbb{Q}).
$$
Note that $N'={\rm dim}\,\mathcal{A}_t=\natural( r^{-1}(t))\,N={\rm
deg}\,(B'/B)\,\,N$. We develop the product above, and observe that
the only nonzero terms appearing in this development come from by
taking $t'=t'''$ and putting the monomial
$$pr_{t'}^*(h_{t'}))^{N-2}\cup_{t''\not=t'}
pr_{t''}^*(h_{t''}))^{N}$$ in the
 middle term. The other terms are $0$ due to the projection formula.
 Let us explain this in the case of only two summands $K_1,\,K_2$ of dimension
 $r$ with polarizations $l_1,\,l_2$
 and
 two cycles $Z_i\in CH^2(S\times K_i)$ giving rise to
 a cycle of the form
 $$(Id_S,pr_1)^*Z_1+(Id_S,pr_2)^*Z_2\in CH^2(S\times K_1\times K_2),$$
 where the $pr_i$'s are the projections from
 $S\times K_1\times K_2$ to $K_i$: Then
 we have
$$[{(pr_1,Id_S)^*}^tZ_1]^*\circ (l_1+l_2)^{2r-2}\cup\circ
 [{(Id_S,pr_2)}^*Z_2]^*=0:H^2(S,\mathbb{Q})_{tr}\rightarrow H^2(S,\mathbb{Q})_{tr}
 $$
 by the projection formula and for the same reason
$$[{(pr_1,Id_S)^*}^tZ_1]^*\circ (l_1+l_2)^{2r-2}\cup\circ
 [{(Id_S,pr_1)}^*Z_1]^*$$
 $$=[{(pr_1,Id_S)^*}^tZ_1]^*\circ \sum_{0\leq k\leq r}\binom{2r-2}{k}(l_1^kl_2^{2r-2-k})\cup\circ
 [{(Id_S,pr_1)}^*Z_1]^*$$
 $$\underset{proj.\,formula}{=}\binom{2r-2}{r}{\rm deg}\,(l_2^r)\,\,[Z_1]^*\circ
 l_1^{r-2}\cup\circ [^tZ_1]^*
 :H^2(S,\mathbb{Q})_{tr}\rightarrow H^2(S,\mathbb{Q})_{tr}.
 $$
 We thus get (as the degrees  ${\rm deg}\,(h_{t'}^{N})$ of the polarizations $h_{t'}$ on the abelian varieties $\mathcal{K}_{t'}$ are all equal):
$$[\Gamma_t]^*=M({\rm deg}\,(h_{t'}^{N}))^{{\rm deg}\,(B'/B)-1}(\sum_{t'\in r^{-1}(t)}[^t\mathcal{Z}'_{t'})]^*)\circ h_{t'}^{N-2}\cup\circ [\mathcal{Z}'_{t'}]^* )
\circ \pi_{tr,t}:
$$
$$H^*(\mathcal{S}_t,\mathbb{Q})\rightarrow
H^*(\mathcal{S}_t,\mathbb{Q}),
$$
where $M$ is the multinomial coefficient appearing in front of the
monomial $x_1^{N-2}x_2^N\ldots x_{{\rm deg}\,(B'/B)}^N$ in the
development of $(x_1+\ldots+ x_{{\rm deg}\,(B'/B)})^{N{\rm
deg}\,(B'/B)-2}$.

By (\ref{autrebonne}) and (\ref{labonne}), we conclude that
$$[\Gamma_t]^*=M({\rm deg}\,(h_{t'}^{N}))^{{\rm deg}\,(B'/B)-1}{{\rm deg}\,(B'/B)} \lambda_g\pi_{tr,t}:H^*(\mathcal{S}_t,\mathbb{Q})\rightarrow
H^*(\mathcal{S}_t,\mathbb{Q}),
$$
which proves formula (\ref{eqcohfin}) with $\lambda=M({\rm
deg}\,(h_{t'}^{N}))^{{\rm deg}\,(B'/B)-1}{{\rm deg}\,(B'/B)}
\lambda_g$.
 \cqfd

 {\bf Proof of Theorem \ref{autreforme}.} Consider the cycle
$\mathcal{Z}\in CH^2(\mathcal{S}\times_B\mathcal{A})_\mathbb{Q}$ of
Lemma \ref{le2fin}, and the cycle $\Gamma\in
CH^2(\mathcal{S}\times_B\mathcal{S})_\mathbb{Q}$
\begin{eqnarray}
\label{formulafinfin} \Gamma:=\Delta_{tr}\circ {^t\mathcal{Z}}\circ
c_1(\mathcal{L})^{N'-2}\circ \mathcal{Z}\circ \Delta_{tr}
\end{eqnarray}
where  now $c_1(\mathcal{L})$ is the class of $\mathcal{L}$ in
$CH^1(\mathcal{A})$ and we denote by $c_1(\mathcal{L})^{N-2}$ the
relative self-correspondence of $\mathcal{A}$ given by the cycle
$c_1(\mathcal{L})^{N-2}$ supported on the relative diagonal of
$\mathcal{A}$ over $B$ (it thus induces the intersection product
with $c_1(\mathcal{L})^{N-2}$ on Chow groups). Furthermore  the
composition of correspondences is the relative composition over $B$,
and $\Delta_{tr}$ is the generic transcendental motive, (which is
canonically defined in our case, at least after restricting to a
Zariski open set of $B$) obtained as follows: Choose a $0$-cycle
$o_{\overline{\eta}}$ of degree $1$ on the generic geometric fiber
$\mathcal{S}_{\overline{\eta}}$ and choose a basis $L_1,\ldots, L_k$
of ${\rm Pic}\,\mathcal{S}_{\overline{\eta}}$. We have then the
projector $\Delta_{alg,\overline{\eta}}\in
CH^2(\mathcal{S}_{\overline{\eta}})_\mathbb{Q}$ defined as in
(\ref{gogo}) using the fact that the intersection pairing on the
group of cycles
$<o_{\overline{\eta}},\,\mathcal{S}_{\overline{\eta}},\,L_i>$ is
nondegenerate. The transcendental projector
$\Delta_{tr,\overline{\eta}}$ is defined as
$\Delta_{\mathcal{S}_{\overline{\eta}}}-\Delta_{alg,\overline{\eta}}$.
This is a codimension $2$ cycle on
$\mathcal{S}_{\overline{\eta}}\times_{\overline{\mathbb{C}(B)}}\mathcal{S}_{\overline{\eta}}$
but it comes from a cycle on
$\mathcal{S}''\times_{B''}\mathcal{S}''$ for some generically finite
covers $B''\rightarrow B,\,\mathcal{S}''=\mathcal{S}\times_BB''$,
and the latter descends finally to a codimension $2$-cycle on
$\mathcal{S}\times_B\mathcal{S}$, which one checks to be a multiple
of  a projector, at least over a Zariski open set of $B$.

The cycle
$\Gamma$ satisfies (\ref{eqcohfin}),
which we rewrite as
\begin{eqnarray}
\label{eqfinfin} [\Gamma'_t]=0\,\,{\rm
in}\,\,H^4(\mathcal{S}_t\times \mathcal{S}_t,\mathbb{Q}),
\end{eqnarray}
where $\Gamma':=\Gamma
-\lambda\Delta_{tr}$.

Lemma \ref{le1fin} tells us that the fibered self-product
$\mathcal{S}\times \mathcal{S}$ is rationally connected. We can thus
apply Theorem \ref{theomain} and conclude that $\Gamma'_t\in
CH^2(\mathcal{S}_t\times \mathcal{S}_t)_\mathbb{Q}$ is nilpotent. It
follows that $\Gamma_t=\lambda\Delta_{tr,t}+N_t$, where
$\lambda\not=0$ and $N_t$ is a nilpotent cycle in
$\mathcal{S}_t\times \mathcal{S}_t$  having the property that
$$\Delta_{tr,t}\circ N_t=N_t\circ \Delta_{tr,t}=N_t.$$
 We conclude immediately from the standard inversion formula for $\lambda\,\mathbb{I}+N$, with $N$ nilpotent and $\lambda\not=0$, that there exists a correspondence
 $\Phi_t\in CH^2(\mathcal{S}_t\times \mathcal{S}_t)_\mathbb{Q}$ such that
 $$\Phi_t\circ \Gamma_t=\Delta_{tr,t}\,\,{\rm in}\,\ CH^2(\mathcal{S}_t\times \mathcal{S}_t)_\mathbb{Q}.$$
Recalling   now  formula (\ref{formulafinfin}):
$$\Gamma_t:=\Delta_{tr,t}\circ{^t\mathcal{Z}_t}\circ c_1(\mathcal{L}_t)^{N'-2}\circ \mathcal{Z}_t\circ \Delta_{tr,t},$$
we proved (\ref{formuleaprouver}) with $$Z'=\Phi_t\circ
\Delta_{tr,t}\circ{^t\mathcal{Z}_t}\circ
c_1(\mathcal{L}_t)^{N'-2}.$$

\cqfd
Let us finish by  explaining the following corollaries.
They all follow  from Kimura's theory of finite dimensionality and are a strong
motivation to establish the Kuga-Satake correspondence at a Chow theoretic
level rather than cohomological one.
\begin{coro} \label{corok3vrai}
With the same assumptions as in Theorem \ref{theocond}, for any $\sigma\in B$,
any self-correspondence of $S_\sigma$
which is homologous to $0$ is nilpotent. In particular, for any finite group
action $G$ on $S_\sigma$ and any projector
$\pi\in \mathbb{Q}[G]$, if $H^{2,0}(S_\sigma)^\pi=0$, then
$CH_0(S_\sigma)_{\mathbb{Q},hom}^\pi=0$.

\end{coro}
 {\bf Proof.} The first statement follows from Theorem
 \ref{theocond} Kimura's work
 (cf. \cite{kimura}) which says that abelian varieties have a finite dimensional motive
 and that for any finite dimensional motive, self-correspondences homologous to $0$ are nilpotent.

 In the case of a finite group action, we consider as in the previous section the self-correspondence
 $\Gamma^\pi$. It is nonnecessarily homologous to $0$, but as it acts as $0$ on
 $H^{2,0}(S)$, its class can be written as
 $$[\Gamma^\pi]=\sum_i\alpha_i [C_i]\times [C'_i]\,\,{\rm in}\,\, H^4(S_\sigma\times S_\sigma,\mathbb{Q})$$
 for some rational numbers $\alpha_i$ and curves $C_i,\,C'_j$ on $S_\sigma$.
 Then we have
 $$[\Gamma^\pi-\sum_i\alpha_i C_i\times C'_i]=0\,\,{\rm in}\,\, H^4(S_\sigma\times S_\sigma,\mathbb{Q}),$$
 from which we conclude by Kimura's theorem that the self-correspondence
 $$Z:=\Gamma^\pi-\sum_i\alpha_i C_i\times C'_i\in CH^2(S_\sigma\times S_\sigma)_\mathbb{Q}$$ is nilpotent.
 It follows that $Z_*:CH_0(S_\sigma)_{\mathbb{Q},hom}\rightarrow CH_0(S_\sigma)_{\mathbb{Q},hom}$ is nilpotent. As it is equal to
 $$\Gamma^\pi_*:CH_0(S_\sigma)_{\mathbb{Q},hom}\rightarrow CH_0(S_\sigma)_{\mathbb{Q},hom},$$
 which is  the projector on $CH_0(S_\sigma)_{\mathbb{Q},hom}^\pi$, we conclude that
 $CH_0(S_\sigma)_{\mathbb{Q},hom}^\pi=0$.

 \cqfd
\begin{coro}\label{corofinalvrai} With the same assumptions as in Theorem \ref{theocond}, the transcendental part
of the Chow motive of any member of the family of $K3$ surfaces
parameterized by $\mathbb{P}(H^0(X,E))$ is indecomposable, that is,
any submotive of it is either the whole motive or the $0$-motive.
\end{coro}
{\bf Proof.} Recall that the transcendental motive of $S_\sigma$ is
$S_\sigma$ equipped with the projector $\pi_{tr}$ defined in
(\ref{gogo}). Let now $\pi\in CH^2(S_\sigma\times
S_\Sigma)_\mathbb{Q}$ be a projector of the transcendental motive of
$S_\sigma$, that is $\pi\circ \pi_{tr}=\pi_{tr}\circ\pi=\pi$. Since
$h^{2,0}(S_\sigma)=1$, $\pi_*$ acts either as $0$ or as $Id$ on
$H^{2,0}(S_\sigma)$. In the first case, ${\rm Ker}\,(\pi_*)_{\mid
H^2(S,\mathbb{Q})_{tr}}$ is a sub-Hodge structure with
$(2,0)$-component equal to $H^{2,0}(S)$. Its orthogonal complement
is then contained in $NS(S_\sigma)_\mathbb{Q}$ which implies that
$\pi_*=0$ on $H^2(S,\mathbb{Q})_{tr}$. In the second case, we find
similarly that $\pi_*=Id$ on $H^2(S,\mathbb{Q})_{tr}$. Since
$\pi=\pi_{tr}\circ\pi=\pi\circ \pi_{tr}$,
 it follows that
 $\pi_*$ acts either by $0$ or as $\pi_{tr}$ on
$H^*(S,\mathbb{Q})$. Hence  the cohomology class of either $\pi$ or
$\pi_{tr}-\pi$ is equal to $0$, from which we conclude by theorem
\ref{theocond} that $\pi$ or $\pi_{tr}-\pi$ is nilpotent. As both
are projectors, we find that $\pi=0$ or $\pi_{tr}=0$ in
$CH^2(S_\sigma\times S_\sigma)_\mathbb{Q}$.

\cqfd

Centre de math\'ematiques Laurent Schwartz

91128 Palaiseau C\'edex, France

\smallskip
 voisin@math.jussieu.fr
    \end{document}